\documentclass[12pt]{amsart}
\usepackage{amssymb}
\usepackage{amsfonts}
\usepackage{latexsym}

\input xy
\xyoption{matrix}\xyoption{arrow}\xyoption{curve}

\long\def\ignore#1{#1}

\newtheorem{theorem}{Theorem}[section]
\newtheorem{observation}[theorem]{Observation}

\newtheorem{corollary}[theorem]{Corollary}
\newtheorem{lemma}[theorem]{Lemma}
\newtheorem{proposition}[theorem]{Proposition}
\def\Hom{\mathop{\rm Hom}\nolimits}

\def\Ext{\mathop{\rm Ext}\nolimits}
\def\Tor{\mathop{\rm Tor}\nolimits}

\def\pd{\mathop{\rm p\,dim}\nolimits}
\def\pdim{\mathop{\rm pdim}\nolimits}
\def\Findim{\mathop{\rm Fin.dim}\nolimits}
\def\findim{\mathop{\rm fin.dim}\nolimits}

\def\ker{\mathop{\rm Ker}\nolimits}

\def\coker{\mathop{\rm Coker}\nolimits}
\def\Im{\mathop{\rm Im}\nolimits}

\def\Mod{\mathop{\rm Mod}\nolimits}

\def\add{\mathop{\rm add}\nolimits}
\def\Add{\mathop{\rm Add}\nolimits}

\def\dirlim{\mbox{\,\hbox{lim}\kern-1.5em
                  \lower1.5ex\hbox{$\longrightarrow$}\,}}
\def\nat{{\Bbb N}}
\def\qed{{\ $\square$}}
\def\lamod{\mbox{$\Lambda$-{\rm mod}}}
\def\laMod{\mbox{$\Lambda$-{\rm Mod}}}
\def\pinf{\mbox{${\mathcal P}^\infty(\Lambda$-{\rm mod}$)$}}
\def\Pinf{\mbox{${\mathcal P}^\infty(\Lambda$-{\rm Mod}$)$}}
\def\lim{\mathop{\rm lim}\nolimits}

\def\lFindim{\mathop{\rm \ell Fin.dim}\nolimits}
\def\lfindim{\mathop{\rm \ell fin.dim}\nolimits}
\def\modla{\mbox{{\rm mod}-$\Lambda$}}
\def\pinfrt{\mbox{${\mathcal P}^\infty(\modla)$}}
\def\filt{\mathop{\rm filt}\nolimits}
\def\pdim{\mathop{\rm p\,dim}\nolimits}

\begin{document}

\title[Finite and Infinite Dimensional Representations of Algebras]{A Homological Bridge Between Finite and Infinite Dimensional Representations of Algebras}

\author[B. Huisgen-Zimmermann]{B. Huisgen-Zimmermann}
\address{
Department of Mathematics \\
University of California\\
Santa Barbara, CA 93106 \\
U.S.A.
}
\email{birge@math.ucsb.edu}
\author[S. O. Smal\o]{S. O. Smal\o}
\address{
Department of Mathematics \\
The Norwegian University for Science and
Technology  \\
7055 Dragvoll \\
Norway
}
\email{sverresm@math.ntnu.no}

\thanks{The work of the first author
was partially supported by an NSF grant, while that of the second author was
partially supported by a Fulbright grant.}

\begin{abstract}
Given a finite dimensional algebra $\Lambda$, we show that a
frequently  satisfied finiteness condition for the category $\pinf$
of all finitely generated (left) $\Lambda$-modules of finite projective
dimension, namely contravariant finiteness of $\pinf$ in $\lamod$, forces
arbitrary modules of finite projective dimension to be direct limits of
objects in
$\pinf$.  Among numerous applications, this yields an encompassing sufficient
condition for the validity of the first finitistic dimension conjecture,
that is, for the little finitistic dimension of
$\Lambda$ to coincide with the big (this is well-known to fail over finite
dimensional algebras in general). 
\end{abstract}

\maketitle

\setcounter{section}{1} \setcounter{theorem}{0}
\begin{center}{\bf 1. Introduction}\end{center}
\medskip

Let $\Lambda$ be an Artin algebra -- the reader may, for example, think of a
finite dimensional algebra over a field -- and let $\lamod$ be
the category of all {\it finitely generated} left
$\Lambda$-modules, $\laMod$ the category of {\it all} left $\Lambda$-modules.
The full subcategories consisting of the objects of finite projective dimension
in
$\lamod$ and
$\laMod$ are denoted by
$\pinf$ and $\Pinf$, respectively. Our interest is focused on the following
two homological dimensions of $\Lambda$: the left little finitistic dimension,
$\ell\findim\Lambda$, which is the supremum of the projective dimensions
attained on the objects in
$\pinf$, and the left big finitistic dimension, $\ell\Findim\Lambda$, which is
defined correspondingly based on $\Pinf$. For several decades, the smoothest
possible connection between
$\pinf$ and $\Pinf$, namely equality of the little and big finitistic
dimensions, had been conjectured to hold, at least when $\Lambda$ is finite
dimensional over a field. However, in 1991 this equality was shown to fail,
even over finite dimensional monomial relation algebras \cite{HZ:2}. Our goal
here is to prove that a certain, very frequently satisfied, finiteness
condition on
$\pinf$ entails that each object in the big category $\Pinf$ is a direct limit
of objects in the small category $\pinf$; in particular, this condition implies
that the two finitistic dimensions coincide.

Before we state this result with greater precision, we give a brief
overview over previously established connections between $\lamod$ and
$\laMod$. In the late seventies, Ringel used the thorough understanding of 
$\lamod$ for
tame hereditary algebras $\Lambda$ which was available at that point 
to give detailed descriptions
of representations in $\laMod$ in this situation \cite{Rin}. We quote from his
introduction: ``The recent progress in the representation theory of finite
dimensional algebras was limited mainly to the modules of finite length and
one would be interested to know in which way the structure of the modules of
finite length determines the behavior of arbitrary modules. Two results of
this type are known...''. The results he subsequently refers to are as
follows:  If $\Lambda$ has finite
representation type, then all objects in $\laMod$ are direct sums of objects
from $\lamod$ (see \cite{Aus:1} and \cite{RT}); if, on the other hand,
$\Lambda$ fails to have finite representation type, $\laMod\setminus\lamod$
contains indecomposable objects \cite{Aus:2}. In the meantime, Crawley-Boevey
added several striking results to this list by showing that, for a finite
dimensional algebra $\Lambda$ over an algebraically closed base field, certain
non-finitely generated modules, called generic modules, completely
determine the representation type of $\Lambda$; they are the non-finitely
generated indecomposable modules which have finite length over their
endomorphism rings (see \cite{CB:1} and \cite{CB:2}). In particular, he
proved that
$\Lambda$ has infinite representation type if and only if there exist generic
objects in $\laMod$, 
tame representation type if and only if, for each dimension $d$, there
exist only finitely many generic $\Lambda$-modules having 
length $d$ over their endomorphism
rings;  finally, $\Lambda$ has domestic type precisely when there are only
finitely many isomorphism classes of generic modules in $\laMod$ altogether.

In 1979, Auslander and Smal\o\ introduced the concept of contravariant
finiteness of a full subcategory $\mathcal A$ of $\lamod$ \cite{AS}.  $\mathcal A$ is
said to be {\it contravariantly finite} in $\lamod$ in case each object $M$ in
$\lamod$ has an {\it $\mathcal A$-approximation} in the following sense: there
exists a homomorphism $f : A\rightarrow M$ with $A\in\mathcal A$ such that each
map in $\Hom_\Lambda(B,M)$ with $B\in\mathcal A$ factors through $f$. One of the
reasons why this concept is of major importance to the homology of a finite
dimensional algebra lies in the fact that, if $\pinf$ is contravariantly
finite in $\lamod$, then the minimal $\pinf$-approximations
of the simple left $\Lambda$-modules  -- 
$A_1,\dots,A_n$ say --  constitute the
basic building blocks for arbitrary representations in $\pinf$.  In fact, due
to Auslander and Reiten \cite{AR}, a module $M$ belongs to $\pinf$ if and
only if
$M$ is a direct summand of a module $N$ having a filtration $N=N_0\supset
N_1\supset
\dots\supset N_m=0$ such that all consecutive factors $N_i/N_{i+1}$ belong to
$\{A_1,\dots,A_n\}$. In particular, this structure result implies of course that
$\ell\findim\Lambda=
\sup\{\pdim A_i\mid i=1,\dots,n\}$.

Our main result (Theorem 4.4) states that contravariant finiteness of
$\pinf$ in
$\lamod$ forces each left $\Lambda$-module of finite projective dimension to
be a direct limit of finitely generated modules of finite projective
dimension.  More strongly, each object in $\Pinf$ is then a direct limit of
finitely generated modules having filtrations with consecutive factors in
$\{A_1, \dots , A_n\}$.  As a consequence, contravariant
finiteness of
$\pinf$ entails the equality
$\ell\findim\Lambda=
\ell\Findim\Lambda$.  This answers in the positive a question left open in
\cite{HHZ}.  As an extra
bonus, the result tells us where to look if we wish to determine the big
finitistic dimension, a notorious quandary in general. In view of the
preceding paragraph, we can obviously conclude that
$\ell\Findim\Lambda=
\sup\{\pdim A_i\mid i=1,\dots,n\}$ in the contravariantly finite case.  We 
illustrate a `typical' computation of this ilk with an example. Another
application of our main theorem, combining it with a result of
Crawley-Boevey, was pointed out to us by H.~Krause. Namely, contravariant
finiteness of $\pinf$ in $\lamod$ implies covariant finiteness and thus
secures the existence of almost split sequences in $\pinf$. On the other hand,
the conclusion of our result is not left-right symmetric. In the final section
of this paper, we present a finite dimensional monomial relation algebra
$\Lambda$ such that
${\mathcal P}^{\infty}$ is contravariantly finite on one side, while the big and
little finitistic dimensions of $\Lambda$ differ on the other.

Some of our tools should be of independent interest. In Section 2, we show
that, over a left noetherian ring, each left module $M$ is the directed union
of those countably generated submodules which have projective dimensions
bounded above by that of $M$. As a consequence, the left big finitistic
dimension of a left artinian ring equals the supremum of those finite
projective dimensions which are attained on countably generated left modules.
In Section 3, we collect some closure properties of the category
$\Add M$ consisting of all direct sums of direct summands of a
fixed module $M$. In particular, we explore closure of such categories under
direct limits.

{\it Acknowledgements.\/} This paper was written while the second author was
visiting the University of California at Santa Barbara, and he would like to 
thank the members of the mathematics department of UCSB, and in particular his
coauthor, for the hospitality he received during his stay there.  Both
authors are indebted to Hans Sverre Smal\o\ for an excellent job of
typesetting and graphing. 

\bigskip
\setcounter{section}{2} \setcounter{theorem}{0}
\begin{center} {\bf 2. The big finitistic dimension equals the countably 
generated finitistic dimension}\end{center}
\medskip

The primary goal of this section is to show that, for each left artinian ring
$R$, the  big finitistic  dimension, $\ell \Findim R$, coincides with the 
supremum of 
those finite projective dimensions which are attained on countably 
generated left $R$-modules. The argument is an off-spring of Kaplansky's
classical theorem \cite{Kap} stating that every projective R-module is a
direct  sum of countably generated components. 

Recall that a module $X$ is said to be countably presented if there exists 
an exact sequence
$$R^{(K)} \rightarrow R^{(L)} \rightarrow X \rightarrow 0$$
with $|K|$, $|L| \leq \aleph_0$. Clearly, every left noetherian ring has 
the property that each of its countably generated left modules is 
countably presented; the same is true if the base ring $R$ is countable 
or a countably generated algebra over a field.

\begin{proposition}
Let $R$ be a ring with the property that each countably generated left 
$R$-module is countably presented. Then each left $R$-module $M$ is the 
directed 
union of those countably generated submodules which have projective 
dimensions bounded above by that of $M$.\end{proposition}

\noindent {\it Remarks.} 1. This proposition is akin to Corollary 3.2.5 of
\cite{RG}, but gives more structural information about the module $M$. We
thank H.~Krause for bringing this reference to our attention.

2. Observe that, in the conclusion of Proposition 2.1, the attribute
`countably  generated' cannot be replaced by `finitely generated'. In fact,
not even over a finite dimensional algebra, need modules of finite
projective dimension be direct limits of {\it finitely generated} modules
of finite projective dimension. This is an immediate consequence of the 
fact that the big finitistic dimension may exceed the little.
\medskip

{\em Proof of the proposition:}\ \ 
We may assume that $\pdim M < \infty$, say $\pdim M = m$. It clearly 
suffices to show that each countable subset of $M$ 
is contained in a countably generated 
submodule  $M'$ with $\pdim M' \leq m$. 
Let 

\ignore{
$$\xy\xymatrix{
 0\ar[r]& P_m\ar[r]^-{f_m}& P_{m-1}\ar[r]^-{f_{m-1}}&
 \cdots\ar[r]& P_1\ar[r]^-{f_1}& P_0\ar[r]^-{f_0}& M\ar[r]& 0
}\endxy $$
}

\noindent be a projective 
resolution of $M$. By Kaplansky's theorem \cite{Kap}, each $P_i$ is of the
form 
$P_i = \bigoplus_{j\in A_i}P_{ij}$, where all of the $P_{ij}$'s are countably 
generated.

Given a countable subset $U$ of $M$, choose a countable subset 
$B_0^{(1)}\subseteq A_0$ such that $f_0(
{ \bigoplus_{j\in{B_0^{(1)}}}} P_{0j}) $ 
contains 
$U$. Then the kernel of the restriction  of $f_0$ to $\bigoplus_{j\in 
B_0^{(1)}} P_{0j}$ 
is in turn countably generated by hypothesis, which permits us to pick a
countable set $B_1^{(1)} \subseteq A_1$, with the property  that 
$$\ker(f_0) \cap \biggl( \bigoplus_{j\in{B_0^{(1)}}} P_{0j}\biggr)  
 \subseteq f_1 \biggl(\bigoplus_{j\in B_1^{(1)}} P_{ij}\biggr).$$
An obvious induction further yields countable subsets $B_i^{(1)} \subseteq 
A_i$ 
for $1 \leq i \leq m$  such that    
$$\ker(f_{i-1}) \cap \biggl( \bigoplus_{j\in 
B_{i-1}^{(1)}} 
P_{i-1,j} \biggr)  \subseteq f_i \biggl(\bigoplus_{j\in B_i^{(1)}}
P_{ij}\biggr)$$  
for 
$ 1 \leq i \leq m$.

Now set $B_m^{(2)} = B_m^{(1)}$. In view of the fact that the image 
$$f_m \biggl(\bigoplus_{j\in B_m^{(2)}} P_{mj}\biggr) \subseteq  \ker
f_{m-1}$$  
is countably generated, we can find 
a countable subset $B_{m-1}^{(2)} \subseteq A_{m-1}$ containing 
$B_{m-1}^{(1)}$ 
such that 
$$f_m\biggl(\bigoplus_{j\in B_m^{(2)}} P_{mj}\biggr)  \subseteq  \ker
(f_{m-1}) \cap \biggl(
\bigoplus_{j\in B_{m-1}^{(2)}} P_{m-1,j}\biggr).$$  
 An induction
analogous to  the preceding 
one then yields countable subsets $B_i^{(2)} \subseteq A_i$,  $0 \leq i \leq
m$  such that $B_i^{(1)} \subseteq B_i^{(2)}$ and 
$$f_{i+1}\biggl(\bigoplus_{j\in B_{i+1}^{(2)}} 
P_{ij}\biggr) \subseteq  \ker(f_i) \cap \biggl( \bigoplus_{j\in B_i^{(2)}}
P_{ij}\biggr)$$ 
for 
$0\leq i\leq m-1$.

Next set $B_0^{(3)} = B_0^{(2)}$, and continue. We iterate these 
$m$-step inductions moving back and forth along the given projective 
resolution of $M$, and an induction on this level will supply us, for each 
$  k\in {\bf N}$, with countable sets $B_0^{(k)}, \cdots , B_m^{(k)}$ 
such that    
$ B_i^{(k-1)} \subseteq B_i^{(k)} \subseteq A_i$ for $0 \leq i \leq m$, having
the additional properties that
$$\ker(f_{i-1}) \cap \biggl( \bigoplus_{j\in
B_{i-1}^{(k)}} P_{i-1,j}\biggr) \subseteq  f_i\biggl(\bigoplus_{j\in
B_i^{(k)}}P_{ij}\biggr)$$  
for $1\leq i\leq m$ and odd $k$, as well
as
$$f_{i+1}\biggl(\bigoplus_{j\in B_{i+1}^{(k)}}P_{i+1,j}\biggr)
 \subseteq 
\ker(f_{i}) \cap \biggl( \bigoplus_{j\in B_{i}^{(k)}} P_{i,j}\biggr)$$ 
for $0
\leq i\leq m-1$ and even $k$. 

Finally, we set $B_i=\cup_{k\in {\bf N}}B_i^{(k)}\subseteq A_i$ for $0\leq 
i\leq m$, define $M'$ to be the countably generated submodule 
$f_0(\bigoplus_{j\in B_0}P_{0j})$ of $M$, and observe that 

\ignore{
$$\xy\xymatrix{ 
0\ar[r] &
\bigoplus_{j\in B_m}P_{mj}\ar[r]^-{f_m}
&\cdots\ar[r]&
\bigoplus_{j\in B_1}P_{1j}\ar[r]^-{f_1}
&\bigoplus_{j\in B_0}P_{0j}\ar[r]^-{f_0}
&M'\ar[r]&0}\endxy$$
}

\noindent is a projective resolution of $M'$. This guarantees that $\pdim
M'\leq m$  as required. \qed

\begin{corollary}
If $R$ is a left artinian ring, then 
$\ell \Findim R$ equals 
$$\sup\{\pdim M\mid M {\rm\ a\
countably\ generated\ left\ } R{\rm -module\ with\ } \pdim M<\infty\}.$$
\end{corollary}

{\em Proof:\ }
Let $N$ be any left $R$-module of finite projective dimension. Since $R$ 
satisfies the hypothesis of Proposition 2.1, we have $N= 
\dirlim_{i \in I} N_i$, where $(N_i)_{i \in I}$ is a directed 
family of 
countably generated modules with $\pdim N_i \leq \pdim N$ for $i \in I$.
In view of the fact that the functor $\Tor$ commutes with direct limits, this
implies that the flat  dimension of $N$ is bounded above by the supremum of
the flat dimensions  of the $N_i$. But our hypothesis on $R$ makes the flat
dimension equal to  the projective dimension and thus completes the argument.
\qed

\bigskip
\setcounter{section}{3} \setcounter{theorem}{0}
\begin{center} {\bf 3. Remarks on the category $\Add(M)$}\end{center}
\medskip

Given a module  $M$, we denote by $\Add(M)$ the full subcategory of  
$R$-$\Mod$
having as objects the modules which are isomorphic to arbitrary direct
sums of direct summands of $M$. The main purpose of this section is to 
assemble a number of properties of $\Add(M)$ for a $\Sigma$-pure injective
module
$M$, to be used towards our main theorem. (That $M$ be $\Sigma$-pure
injective means that all direct sums of copies of $M$ are pure injective.)  
However, due to the fact that this category holds interest in its own right,
we will go a little beyond the requirements of the following section. 

Start by
recalling that every finitely generated module over an Artin algebra is
$\Sigma$-pure injective, since of finite length over its endomorphism
ring.  On the other hand, $\Sigma$-pure injectives  are far from being finitely
generated over the base ring, in general;  in particular, any generic module
has this property.  More generally, recall that the results in \cite{Zim} and
\cite{GJ} give an equivalent characterization of
$\Sigma$-pure  injectivity of an $R$-module $M$ in terms of a descending chain 
condition for certain submodules of $M$ over its endomorphism ring, a condition
which is obviously inherited by pure $R$-submodules of $M$; as a
consequence,
$\Sigma$-pure injectivity is passed on to pure submodules. 

Our first observations address closure of the category $\Add(M)$
under direct limits (we follow the convention of  reserving the term `direct
limit' for colimits over directed index sets).  A twin sibling of the
following fact was proved by Lenzing in \cite{Len}, where it is shown that,
for  any finitely presented module $M$ of finite length, the category
$\Add(M)$ has the  mentioned closure 
property.

\begin{observation}  Let  $R$ be any ring and  $M$ a $\Sigma$-pure
injective left $R$-module.  Then every directed system of objects in $\Add(M)$ 
attains its
direct limit in $\Add(M)$. \end{observation}

{\em Proof:\ }  Let $(A_i, f_{ij})_{i,j \in I, i\le j}$ be any directed 
system
of objects $A_i$ in $\Add(M)$, and denote by $A$ its direct limit. Moreover,
let $\iota_j: A_j \rightarrow \bigoplus_{i \in I} A_i$ be the canonical
embeddings.  It is well known that the short exact sequence
$$0 \rightarrow \sum_{i,j \in I, i<j}(\iota_j f_{ij} - \iota_i)(A_i)
 \rightarrow \bigoplus_{i \in I} A_i \rightarrow  A \rightarrow 0$$
is pure.  Its middle term being $\Sigma$-pure injective by hypothesis, the
sequence therefore splits.  Now each $\Sigma$-pure injective module is a direct
summand of components, each of which has a local endomorphism ring
\cite{Zim}, p.~1100, and consequently the Krull-Remak-Schmidt-Azumaya Theorem
guarantees that all terms of the above exact sequence in turn belong to
$\Add(M)$. \qed

\medskip

Interestingly, the question whether $\Add(M)$ is closed under
direct limits can be decided by testing solely well-ordered chains in general. 
We record this fact as

\begin{observation} Given any left $R$-module $M$, the category
$\Add(M)$ is closed under arbitrary direct limits, provided that it is
closed under direct limits  over well-ordered index sets.
\end{observation}

{\em Proof:\ } Again we let $(A_i, f_{ij})_{i,j \in I, i \leq j}$ be a
directed system of objects in $\Add(M)$ and denote its direct limit by $A$.  
If the index set
$I$ is finite, our claim is trivial, so let us suppose that
$|I| = \aleph_{\alpha}$ for some ordinal number $\alpha$.   We will
prove the observation by a transfinite induction on $\alpha$.  In case 
$\alpha = 0$, we may clearly assume that $I = \{ i_1, i_2, i_3, \dots \}$ and
that $i_j < i_k$ in the partial order of $I$ whenever $j < k$.   Our
hypothesis therefore guarantees that $A$ is in turn an object of $\Add(M)$. 

Now we suppose that $\alpha > 0$ and that for all ordinal numbers $\beta <
\alpha$ our claim is true; in other words, we assume that all direct limits
of directed systems of objects in $\Add(M)$ extending over index sets of 
cardinalities less than $\aleph_{\alpha}$ belong to $\Add(M)$.  We 
write
$\Omega_{\alpha}$ for the first ordinal number of cardinality
$\aleph_{\alpha}$ and index the set $I$ by the ordinals below
$\Omega_{\alpha}$, i.e., we write $I = \{ i_\beta \mid  \beta < \Omega_\alpha
\}$. It will be convenient to identify an ordinal number $\gamma$ with 
the set of all ordinal
numbers strictly smaller than $\gamma$.  Our aim is to show that
$A$ is the direct limit of a directed system
$(C_{\beta}, g_{\beta, \gamma})_{\beta, \gamma < \Omega_{\alpha}, \beta \le
\gamma}$ with the property that all of the objects $C_{\beta}$ belong to
$\Add(M)$.  Since the new index set,
$\Omega_{\alpha}$, is totally ordered, we can then again invoke the hypothesis
to complete the proof.

For that purpose, we observe that, given any infinite subset $I' 
\subseteq I$,
there exists a {\it directed} subset $I''$ of $I$ which contains $I'$ and 
has the
same cardinality as $I'$.  This allows us to express $I$ as the directed union
of a chain of subsets $I_{\gamma}$, $\gamma < \Omega_{\alpha}$, such 
that $I_{\gamma}
\subseteq I_{\delta}$ whenever $\gamma < \delta$, that  $I_{\gamma}$ is 
finite whenever $\gamma$
is finite, and $I_{\gamma}$ has cardinality less than or equal to that of
$\gamma$  whenever $\gamma$ is an
infinite ordinal number.  Indeed, set $I_0 =
\{i_0\}$, let $\delta > 0$, and suppose that, for all $\gamma < \delta$, 
directed
subsets $I_{\gamma}$ satisfying the above requirements and having the 
property that $\{
i_\beta \mid \beta < \gamma\} \subseteq I_{\gamma}$ have
already been constructed.  Set
$L_{\delta} = \cup_{\gamma < \delta} I_{\gamma}\cup\{i_{\delta}\}$.  By
construction $L_{\delta}$ is either finite or else has a cardinality 
bounded above by that of
$\delta$.  If
$L_{\delta}$ is finite, we can clearly choose a finite directed subset
$I_{\delta}$ of $I$ containing $L_{\delta}$; if, on the other hand,
$L_{\delta}$ is infinite, the initial remark of this paragraph permits us 
to choose a directed
subset $I_{\delta}$ of $I$ which contains the set $L_{\delta}$ and has 
the same cardinality as the latter.  A subsidiary induction thus gives us 
the desired family
of subsets $I_{\gamma}$ for $\gamma < \Omega_{\alpha}$.  The hypothesis 
of the
principal induction now yields that for each $\gamma$, the direct limit
$C_{\gamma}$ of the directed subsystem $(A_i, f_{ij})_{i,j \in 
I_{\gamma}, i \le j}$ of our
original system is an object of $\Add(M)$.  Letting $g_{\beta,\gamma}: 
C_{\beta} \rightarrow
C_{\gamma}$ for $\beta < \gamma$ be the natural map resulting from the fact
that $I_{\beta} \subseteq I_{\gamma}$, we are thus in a position to apply the
hypothesis to the system $(C_{\beta}, g_{\beta, \gamma})_{\beta, \gamma < 
\Omega_{\alpha},
\beta \le \gamma}$ to conclude that $A = \dirlim C_{\beta}$ in turn 
belongs to
$\Add(M)$. \qed 

\medskip

An additional remark on the category $\Add(M)$ for a
$\Sigma$-pure injective module $M$ will turn out helpful in Section 4, namely

\begin{observation} Let $M$ be $\Sigma$-pure injective, and 
suppose that all
monomorphisms in $\add(M)$ split.  Then the same is true for 
monomorphisms in $\Add(M)$; 
in other words, given any monomorphism $f:A \rightarrow B$ with $A,B$ in 
$\Add(M)$, there
exists a submodule $C$ of $B$ such that $B= \Im(f) \bigoplus C$, and both 
$\Im(f)$ and $C$ again
belong to $\Add(M)$. \end{observation}

{\em Proof:\ } Let $f$ be as specified in the claim, and write $B 
=\bigoplus_{i \in I} B_i$, where
each $B_i$ is isomorphic to a direct summand of $M$.  Then $f$ is 
necessarily pure, since the
image $f(F)$ of any finitely generated submodule $F$ of
$A$ is contained in a finite direct sum of the $B_i$ and hence is a 
direct summand of  $B$ by
hypothesis. As mentioned in the beginning of this section,
$\Sigma$-pure injectivity is passed on to pure submodules, and therefore 
$\Im(f)$ is a direct
summand of $B$.  Finally we use once more the fact that each
$\Sigma$-pure injective module is a direct sum of submodules with local 
endomorphism rings, in
order to deduce that both $\Im(f)$ and any complement of $\Im(f)$ in $B$ are
in  turn direct sums of
direct summands of $M$. \qed

\bigskip
\setcounter{section}{4} \setcounter{theorem}{0}
\begin{center} {\bf 4. $\Findim$ versus $\findim$ in case $\pinf$ is
contravariantly finite}\end{center}
\medskip

This section is devoted to proving our main result, namely that, for any
Artin  algebra $\Lambda$, the left big and little finitistic dimensions
coincide, provided that the full subcategory  $\pinf$ of finitely
generated left $\Lambda$-modules of finite projective dimension is 
contravariantly finite in the category of all finitely generated
left $\Lambda$-modules.

Recall that $\pinf$ is contravariantly finite in $\lamod$ if, for each finitely
generated left $\Lambda$-module $M$, there exists an object $A \in \pinf$ and
a homomorphism $f: A \rightarrow M$ such that the induced sequence of
functors from $\pinf$ to the category of abelian groups,
$$\Hom_{\Lambda}(-, A)|_{\pinf} \rightarrow \Hom_{\Lambda}(-, M)|_{\pinf}
\rightarrow 0,$$
is exact; see \cite{AS}.  In that case, $A$ is called a (right)
$\pinf$-approxiamtion of $M$.  It is well-known that, existence provided, the
$\pinf$-approximations of minimal length of a given module $M \in \lamod$ are
all isomorphic, and hence it makes sense to refer
to {\it the} minimal
$\pinf$-approximation of $M$ in that case. 

It is known that for any left serial algebra $\Lambda$, the category $\pinf$
is contravariantly finite in $\lamod$ \cite{BHZ}; moreover, the minimal
$\pinf$-approximations of the simple left $\Lambda$-modules can be explicitly
described in that case. Initial criteria for contravariant
finiteness in more general situations were developed in
\cite{HHZ}.

Throughout this section, we will abbreviate $\ell\findim\Lambda$ by
$\findim\Lambda$ and $\ell\Findim\Lambda$ by
$\Findim\Lambda$.

\begin{theorem} Suppose that $\pinf$ is contravariantly finite
in $\lamod$. Then 
$$\findim \Lambda = \Findim
\Lambda.$$
\end{theorem}

{\em Proof:\ } Since $\Findim \Lambda$ is the supremum of those
projective dimensions which are attained on countably generated left
$\Lambda$-modules of finite projective dimension by Corollary 2.2, it is
enough to focus on a countably generated module $M$ with $\pd  M <\infty$, say
$\pd M = n$.  Let $M_1\subseteq M_2\subseteq M_3 \subseteq \cdots 
\subseteq M_r \subseteq \cdots 
\subseteq M$ be a chain of  finitely generated submodules such that
$M=\cup_{i \geq 1} M_i$.  For each i, we fix the beginning of a finitely
generated projective resolution of
$M_i$, say
$$0 \rightarrow \Omega_{n,i} \rightarrow P_{n-1,i}\rightarrow 
\cdots\rightarrow 
P_{1,i} \rightarrow P_{0,i} \rightarrow M_i \rightarrow 0$$ and,
over each of the inclusions  $M_i \rightarrow M_{i+1}$,  we choose a
chain morphism consisting of maps 
$f_{k,i,i+1}: P_{k,i} \rightarrow P_{k,i+1}$ for $0 \le k \le n-1$, and 
$f_{n,i,i+1}: \Omega_{n,i} \rightarrow \Omega_{n,i+1}$. Defining
$f_{k,i,j}=f_{k,j-1,j}\circ  \cdots \circ  f_{k,i,i+1}$
for all $j>i$ and $ 0 \leq k \leq n$, we thus obtain a directed system of 
exact 
sequences indexed by the natural numbers. Passing to the direct limit of this
system gives us  a projective resolution of the module $M$ which we label 
$$0\rightarrow\dirlim \Omega_{n,i}\rightarrow P_{n-1}\rightarrow \cdots 
\rightarrow P_1\rightarrow P_0\rightarrow M\rightarrow 0$$ 
In particular, $P_n := \dirlim \Omega_{n,i}$ 
is projective due to the fact that $\pdim M =n$. Clearly, all $P_k$ are
countably generated since our directed system extends over $\nat$.  Therefore,
we can decompose
$P_n$ in the form
$P_n =\bigoplus_{j \in \nat} Q_{n,j}$, where all the $Q_{n,j}$ are finitely
generated.  Denoting by $g_i:\Omega_{n,i} \rightarrow P_n$  
the canonical morphisms, we obtain an increasing sequence $i_1, i_2,
i_3, \dots$  of natural numbers such that $\bigoplus_{j \leq m} Q_{n,j}
\subseteq
\Im g_{i_m}$ for $m \in \nat$.  It is clearly harmless to pass to a suitable
cofinal subsystem of the $\Omega_{n,i}$, which permits us to assume that
$i_k = k$ for all $k$.  Next we observe that each $\Omega_{n,i}$ contains a 
submodule $P_{n,i} = \bigoplus_{j \leq i} P_{n,j,i}$ isomorphic to the
finite direct sum $\bigoplus_{j \leq i} Q_{n,j}$ with the property that
$g_i$ restricts to an isomorphism
$P_{n,j,i} \rightarrow Q_{n,j}$ for all $j \leq i$;  just keep in mind that
the restriction of the map $g_i$ to the preimage of $\bigoplus_{j \leq
i}Q_{n,j}$ splits.  Since $g_{i+1} f_{n,i,i+1}\bigl(P_{n,i}\bigr) =
g_i\bigl(P_{n,i}\bigr)$, the map $f_{n,i,i+1}$ induces a split monomorphism
$P_{n,i} \rightarrow g_{i+1}^{-1}\bigl(\bigoplus_{j \leq i+1}
Q_{n,j}\bigr)$.  Consequently, an obvious induction on
$i$ allows us to choose the $P_{n,j,i}$ in such a way that the squares

\ignore{
$$\xy\xymatrix{
P_{n,i+1} \ar[r]^-{g_{i+1}} &\bigoplus_{j\le i+1} Q_{nj}\\
P_{n,i} \ar[u]^-{f_{n,i,i+1}} \ar[r]^-{g_i}_\cong &\bigoplus_{j\le i} Q_{nj}
\ar@{_(->}[u]
}\endxy$$
}

\noindent commute. In other words, this process yields a directed subsystem

\ignore{
$$\xy\xymatrix{ 
\vdots& \vdots\\ 
P_{n,i+1} \ar[u] \ar@{^(->}[r] &\Omega_{{n,i+1}} \ar[u]\\
P_{n,i}\ar@{^(->}[u] \ar@{^(->}[r] &\Omega_{n,i} \ar[u]\\
\vdots \ar[u] &\vdots \ar[u]
}\endxy$$
}

\noindent of the system $\bigl( \Omega_{n,i}, f_{n,i,j}\bigr)_{i,j \in \nat,
i\leq j}$ such that $\dirlim P_{n,i} = \dirlim \Omega_{n,i} = P_n$.

At this point, we interrupt the argument with a lemma which will allow us to
supplement our original directed system of resolutions 
$${\bf (S_i)}\ \ \ \ \ P_{n-1,i} \rightarrow P_{n-2,i} \rightarrow
\dots \rightarrow P_{0,i} \rightarrow M_i \rightarrow 0$$ by a system
$${\bf (T_i)}\ \ \ \ \  P_{n-1,i}'
\rightarrow P_{n-2,i}' \rightarrow \dots \rightarrow P_{0,i}' \rightarrow N_i
\rightarrow 0,$$ 
where $N_i\in \pinf$, together with an epimorphism ${\bf (T_i)} \rightarrow
{\bf (S_i)}$ of directed systems, to the effect that
$N =\dirlim N_i$ has a projective dimension bounded above by $\findim \Lambda$
and the kernel of the induced epimorphism $N \rightarrow M$ is `under
control'.  

The following lemma is based on the fact that contravariant
finiteness of $\pinf$ in $\lamod$ implies the
existence of an injective cogenerator inside the category $\pinf$ (see
\cite{AS}); indeed, if $I$ is the minimal right
$\pinf$-approximation of the minimal injective
cogenerator for $\lamod$, then each object of $\pinf$ embeds into an object of $\add I$ and every inclusion 
$I' \hookrightarrow X$ with $I'$ in $\add I$ and $X$ in $\pinf$ splits.  In
particular, this entails that every monomorphism in $\add(I)$ splits, i.e.,
the hypotheses of Observation 3.3 are satisfied for $M=I$.
 
\begin{lemma} Let $I$ be a relative injective cogenerator for $\pinf$ as above, and suppose that we are given
an exact commutative diagram of the form

\ignore{
$$\xy\xymatrix{ 
0\ar[rr]&&X\ar[rr]^{\beta}&&P\ar[rr]^{\alpha}
&&Y\ar[rr]&&0\\
&Q\ar[ur]^g\\
0\ar@{-}[r]&\ar[r]&X'\ar[rr]^{\beta'}
\ar[uu]_{f_X}
&&P'
\ar[rr]^{\alpha'}
\ar[uu]_{f_P}
&&Y'\ar[rr]\ar[uu]_{f_Y}
&&0\\
&Q'\ar[uu]^(.7){h_Q}
\ar[ur]_{g'}
}\endxy$$
}

\noindent with $Q$, $Q'$, $P$ and $P' \in \pinf$. Then
there  exist modules $I_0$ and $I_0'$ 
in $\add I$, together with homomorphisms $\gamma:Q \rightarrow I_0$ and $
\gamma' :Q' \rightarrow  I_0'$, as well as a homomorphism    
$h:I_0' \rightarrow I_0$ such that the following diagram has exact rows
and commutes.

\ignore{
$$\xy\xymatrix{ 
&0\ar[rr]&&X\ar[rr]^{\beta}&&P\ar[rr]^{\alpha}&&Y\ar[rr]&&0\\
0\ar[rr]&&Q\ar[ur]^g\ar[rr]^(.6){\left({{\beta 
g}\atop{\gamma}}\right)}&\ar[u]&P\bigoplus 
I_0\ar[ur]^{(1,0)}\ar[rr]&\ar[u]&Z\ar[ur]\ar[rr]&\ar[u]&0\\
 &0\ar@{-}[r] &\ar[r] &X'\ar@{-}[u]^{f_X} \ar@{-}[r] &\ar[r]^{\beta'}&
P'\ar@{-}[u]_{f_P}\ar@{-}[r]&\ar[r]^{\alpha'}&
Y\ar@{-}[u]_{f_Y}\ar[rr]&&0\\
0\ar[rr]&&Q'\ar[uu]^(.7){h_Q}\ar[ur]^{g'}
\ar[rr]^{\left({{\beta' g'}\atop{\gamma'}}\right)}&& P'\bigoplus I_0' 
\ar[uu]^(0.7){\left({{f_P}\atop{0}}{{0}\atop{h}}\right)}\ar[ur]^{(1,0)}
\ar[rr]&&
Z'\ar[uu]\ar[ur]\ar[rr]&&0
}\endxy$$
}

\noindent Here $Z= \coker{\left({{\beta g}\atop{\gamma }}\right)}$,
 $Z'= \coker{ \left({{\beta' g'}\atop{\gamma' }}\right)}$ are in $\pinf$, 
and the remaining maps are induced by the cokernels.
\end{lemma}

{\em Proof of the lemma:\ } Since $I$ is a cogenerator for $\pinf$, we can choose
monomorphisms 
$\gamma:Q \rightarrow I_0$ and  $\gamma':Q' \rightarrow I_0'$.  So if we
introduce maps $Q \rightarrow P \oplus I_0$ and  $Q' \rightarrow P \oplus I_0'$
as in the above diagram and denote by $Z$ and $Z'$ their cokernels,
respectively, the two squares in each of the top and bottom planes
commute.  Moreover, $Z,Z'$ in turn belong to $\pinf$ by  the hypothesis on
$Q$ and $Q'$.  Next we use the relative injectivity of $I$ to obtain $h: I_0'
\rightarrow I_0$ such that
$\gamma h_Q = h \gamma'$.  It is now straightforward to check that the entire
diagram commutes, which completes the proof of the
lemma. \qed
 
\medskip

We return to the proof of Theorem 4.1.  At this point, we label the maps in the
projective resolutions of the $M_i$, say $g_{k,i}: P_{k,i} \rightarrow
P_{k-1,i}$.  By applying the lemma, first to the diagrams 

\ignore{
$$\xymatrixcolsep{1.3pc}
\xy\xymatrix{ 
0 \ar[rr] &&P_{n-1,i+1}/\Omega_{n,i+1}
\ar[rr]^-{g_{n-1,i+1}} &&P_{n-2,i+1} \ar[rr]^-{g_{n-2,i+1}}
&&\Im(g_{n-2,i+1}) \ar[rr] &&0\\ 
 &P_{n-1,i+1}/P_{n,i+1} \ar[ur]\\
0 \ar@{-}[r] &\ar[r] &P_{n-1,i}/\Omega_{n,i} \ar[rr]^-{g_{n-1,i}}
\ar[uu]_{f_{n-1,i,i+1}} &&P_{n-2,i} \ar[rr]^-{g_{n-2,i}}
\ar[uu]_{f_{n-2,i,i+1}} &&\Im(g_{n-2,i}) \ar[rr] \ar[uu] &&0\\
 &P_{n-1,i}/P_{n,i} \ar[uu] \ar[ur]
}\endxy$$
}

\noindent for $i \in \nat$, and by then moving inductively along the sequences
$$P_{n-1,i} \rightarrow P_{n-2,i} \rightarrow \dots \rightarrow P_{0,i}
\rightarrow M_i \rightarrow 0,$$
we obtain the following directed system of
diagrams in $\lamod$: 

\ignore{
$$\xy\xymatrix@C-4mm{ & &0 & 0& &0&0&0\\ 
0\ar[r]&\Omega_{n,i}\ar[r]& P_{n-1,i} \ar[u] \ar[r]^{g_{n-1,i}}
&P_{n-2,i} \ar[r]^{g_{n-2,i}} \ar[u]&
\cdots\ar[r] &P_{1,i}\ar[u]\ar[r]&P_{0,i}\ar[u]\ar[r]&
M_i\ar[u]\ar[r] &0\\ 
0\ar[r] &P_{n,i} \ar[r] \ar@{^(->}[u]&
P_{n-1,i}\ar[r]\ar@{=}[u]&P_{n-2,i}\bigoplus
I_{n-2,i}\ar[r]\ar[u]^{(1,0)}& \cdots\ar[r]&P_{1,i}\bigoplus
I_{1,i}\ar[u]^{(1,0)}\ar[r]&P_{0,i}\bigoplus
I_{0,i}\ar[u]^{(1,0)}\ar[r]& N_i \ar[u] \ar[r]&0\\
 &&0 \ar[u] \ar[r] & I_{n-2,i}\ar[r]\ar[u]^{\left({{0}\atop{1}}\right)}&
\cdots\ar[r]& I_{1,i}\ar[u]^{\left({{0}\atop{1}}\right)}\ar[r]&
I_{0,i}\ar[u]^{\left({{0}\atop{1}}\right)}\ar[r]& K_i \ar[u] \ar[r] &0\\
 & & & 0 \ar[u] & &0\ar[u] &0\ar[u] &0\ar[u]
 }\endxy$$
}

\noindent Here the objects $I_{k,i}$ all belong to $\add(I)$. The upper two
horizontal sequences of each of these diagrams are exact by construction,
whereas the induced kernel sequence in the third row will not be exact in
general; in fact, the homology in the term labeled $n-2$ of that sequence is
isomorphic to  
$\Omega_{n,i}/P_{n,i}$.  On the other hand, the direct limit of the inclusions
$P_{n,i} \hookrightarrow \Omega_{n,i}$ is the identity $P_n \rightarrow P_n$,
and hence the snake lemma ensures that, post-limits, we arrive at an {\it
exact} commutative diagram of the form 

\ignore{
$$\xy\xymatrix@C-4mm{ & &0 & 0& &0&0&0\\ 
0\ar[r]&P_{n}\ar[r]&
P_{n-1}\ar[u] \ar[r]&P_{n-2}\ar[r]\ar[u]&
\cdots\ar[r]&P_{1}\ar[u]\ar[r]&P_{0}\ar[u]\ar[r]&
M\ar[u]\ar[r]&0\\ 
0\ar[r]&P_{n,}\ar[r]\ar@{=}[u]&
P_{n-1}\ar[r]\ar@{=}[u]&P_{n-2}\bigoplus
I_{n-2}\ar[r]\ar[u]^{(1,0)}& \cdots\ar[r]&P_{1}\bigoplus
I_{1}\ar[u]^{(1,0)}\ar[r]&P_{0}\bigoplus
I_{0} \ar[u]^{(1,0)} \ar[r] & N \ar[u]\ar[r]&0\\
 &&0\ar[u] \ar[r]& I_{n-2} \ar[r] \ar[u]^{\left({{0}\atop{1}}\right)} &
\cdots \ar[r] & I_{1} \ar[u]^{\left({{0}\atop{1}}\right)}\ar[r] &
I_{0} \ar[u]^{\left({{0}\atop{1}}\right)} \ar[r] & K \ar[u]\ar[r]&0\\ & &
& 0\ar[u]
& &0\ar[u]
&0\ar[u]
&0\ar[u]
 }\endxy$$
}

\noindent Now Observation 3.1 yields $I_j \in \Add (I)$ and, as a consequence,
Observation 3.3 gives us $K \in \Add (I)$.  Since the $N_i$ belong to 
$\pinf$ by construction, the relative cogenerating property of $I$ furthermore
permits us to choose a directed system of embeddings $N_i \hookrightarrow I_i'$
with $I'_i \in \add (I)$, which shows that $N = \dirlim N_i$ in turn embeds
into an object of $\Add (I)$ by 3.1.  In view of 3.3, this implies that the
exact sequence
$0 \rightarrow K
\rightarrow N \rightarrow M  \rightarrow 0$ splits, and therefore $\pd M \leq
\pd N \leq \sup \{\pd N_i \mid i 
\in \nat \} \leq \findim \Lambda$. \qed
\medskip

In the light of Theorem 4.1, the following is an immediate consequence of the
result of Auslander and Reiten quoted in the introduction; just keep in mind that
each
$\pinf$-approximation of a module $M$ contains a minimal one
as a direct summand \cite{AR}.

\begin{corollary}  If $\pinf$ is contravariantly finite in $\lamod$, and $C_1,
\dots , C_n$ are arbitrary $\pinf$-approximations of the simple left
$\Lambda$-modules, then
$$\Findim \Lambda = \sup \{\pd C_i \mid 1 \leq n \}.$$ \qed
\end{corollary}

We can actually strengthen the conclusion of Theorem 4.1, so as to provide  
information on the structure of the non-finitely generated objects in
$\Pinf$ as follows:
Given modules $M_1,\cdots, M_n$ in $\lamod$ we denote by $\filt 
(M_1,\cdots, M_n)$ the full subcategory of $\lamod$ having as 
objects all finitely generated modules $X$ that possess filtrations with
consecutive factors  in $\{M_1,\cdots,M_n\}$; in other words, $\filt(M_1,
\cdots, M_n)$  consists of those modules $X$ which contain chains of the form 
$X=X_0\supseteq X_1 \supseteq \cdots \supseteq X_r =0$ such that each 
factor $X_i/X_{i+1}$ is isomorphic to some $M_j$. Moreover,
$\overrightarrow\filt ( M_1, 
\cdots, M_n)$ will stand for the full subcategory of $\laMod$ the 
objects of which are the direct limits of modules in $\filt (M_1, 
\cdots, M_n)$. Due to the
above-mentioned result of Auslander and 
Reiten, we know: In case $\pinf$ is contravariantly finite in $\lamod$ and 
$A_1, \cdots, A_n$ are the minimal $\pinf$-approximations of the 
simple left $\Lambda$-modules, $\pinf$ consists of the direct summands of 
modules in $\filt (A_1, \cdots, A_n)$. 
In view of the proof of Theorem 4.1,
this description of the finitely 
generated modules of finite projective dimension extends to non-finitely 
generated candidates as follows.

\begin{theorem}
If $\pinf$ is contravariantly finite in $\lamod$ and $A_1,\cdots,A_n$  
are as above, then 
$$\Pinf= \overrightarrow\filt (A_1, \cdots, A_n).$$
In particular, each object of $\Pinf$ is a direct limit of modules in $\pinf$.
\end{theorem}

{\em Proof\,:\ } That $\overrightarrow\filt(A_1,\dots,A_n)$ is contained in
$\Pinf$ is clear. For the other inclusion, start by noting that each full
subcategory
$\mathcal C$ of
$\laMod$ which is closed under direct limits is also closed under direct
summands. This is well known (see e.g. \cite{CB:3}, Lemma 1), but we include
the easy argument: If $B$ is a direct summand of an object $C$ in $\mathcal C$ and
$\pi : C\rightarrow C$ a projection onto $B$, then $B$ is the direct limit of
the system 

\ignore{
$$\xy\xymatrix{
C \ar[r]^\pi &C \ar[r]^\pi &C \ar[r]^\pi &\cdots
}\endxy$$
}

By Proposition 2.1, each
object in
$\Pinf$ is the direct union of countably generated objects in $\Pinf$, whence
it suffices to show that each countably generated module $M$ of finite
projective dimension belongs to
$\overrightarrow \filt (A_1,\dots,A_n)$. But the proof of Theorem 4.1 shows
that $M$ is a direct summand of a direct limit of modules in $\pinf$ and thus
belongs to the closure of ${\pinf}$ under direct limits as explained above.
That the latter category is contained in
$\overrightarrow \filt (A_1,\dots,A_n)$, finally, follows from \cite{AR} and
another application of our initial remark. \qed\medskip 

Note that our arguments for Theorems 4.1 and  4.4 only use the 
existence of a 
$\pinf$
approximation of 
$D(\Lambda)$,
where $D$ is the standard  duality 
$\modla \rightarrow \lamod$.
This observation does not lead to any significant generalization of our 
results however. Indeed, by \cite{HU}, in the presence of the 
inequality 
$\findim \Lambda <\infty$,
the existence of a 
$\pinf$-approximation of
$D(\Lambda)$
forces $\pinf$ to be contravariantly finite.

The final consequence of our main results was observed by Henning Krause, who
pointed out to us that a result of Crawley-Boevey applies to our context.

\begin{corollary} If the subcategory $\pinf$ of $\lamod$ is contravariantly
finite, then it is also covariantly finite.  In particular, $\pinf$ has
almost split sequences in that case.
\end{corollary}

{\em Proof\,:\ } Suppose that $\pinf$ is contravariantly finite in $\lamod$.
As we know, $\lFindim\Lambda$ is finite in that case, equal to $d$ say. In
particular, $\Pinf$ is closed under direct products -- just keep in mind that
direct products of projectives are projective. By Theorem 4.4, this is the
same as to say that the closure of $\pinf$ under direct limits is also closed
under direct products.  Hence Theorem 4.2 of \cite{CB:3} tells us that
$\pinf$ is indeed covariantly finite. \qed\medskip

On the other hand, it is not true in general that covariant finiteness of
$\pinf$ in
$\lamod$ implies contravariant finiteness. Indeed, since the category of
modules of projective dimension $\le1$ is always covariantly finite in
$\lamod$ by \cite{AR}, the examples in \cite{IST} and \cite{Sma}
exhibit covariantly finite categories $\pinf$ which fail to be contravariantly
finite. A further discussion of the two properties can be found
in \cite{HZS}.

\bigskip
\setcounter{section}{5} \setcounter{theorem}{0}
\begin{center} {\bf 5. Examples}\end{center}
\medskip

Throughout this section, $\Lambda=K\Gamma/I$ will be a finite dimensional 
path 
algebra modulo relations over a field $K$. We will denote by $J$ its 
Jacobson radical, by $e_i$ the primitive idempotents of $\Lambda$ going 
with the 
vertices of the quiver  $\Gamma$, and by $S_i=\Lambda e_i/Je_i$ 
representatives of the simple left 
$\Lambda$-modules. Moreover,  given a left $\Lambda$-module $M$ and a 
primitive idempotent $e_i$, we will call an element $ x\in M $ a 
{\it top 
element of type $e_i$} in case $ x \in M \setminus JM  $ and $e_ix=x$. We 
start 
with an example which illustrates the applicability of corollary 4.3.
\medskip

\noindent {\bf Example 5.1.}
Let $\Lambda$ be the monomial relation algebra $ K\Gamma/I$ where 
$\Gamma$ is the quiver

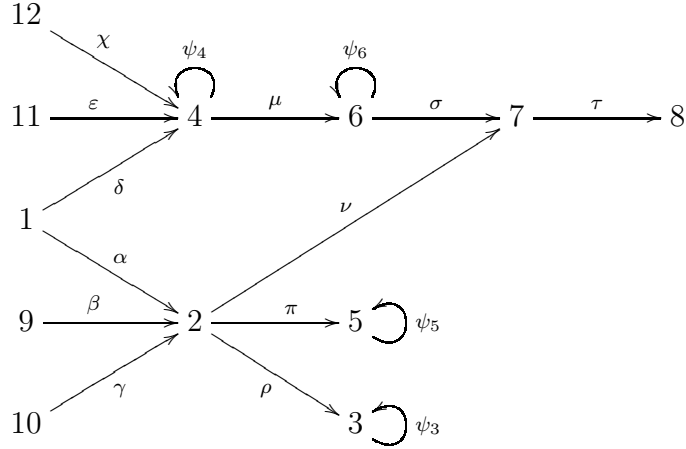
\ignore{
$$\xymatrixcolsep{4pc}
\xy\xymatrix{
12\ar[dr]^(0.4)\chi\\
11\ar[r]^(0.4)\varepsilon&4\ar@(ur,ul)[]_{\psi_4}\ar[r]^\mu&
6\ar@(ur,ul)[]_{\psi_6}\ar[r]^\sigma&7\ar[r]^\tau&8\\
1\ar[ur]_\delta \ar[dr]^\alpha \\
9 \ar[r]^(0.4)\beta&2 \ar[uurr]^\nu \ar[r]^(0.6)\pi \ar[dr]_\rho &5 
\ar@(dr,ur)[]_{\psi_5}\\
10 \ar[ur]_\gamma &&3 \ar@(dr,ur)[]_{\psi_3} 
}\endxy$$
}

\noindent and $I \subseteq K\Gamma$ is the unique ideal of the path
algebra with the property that the  indecomposable projective left
$\Lambda$-modules have the following  graphs. 
For our graphing conventions, we refer the reader to \cite{HZ:1}.

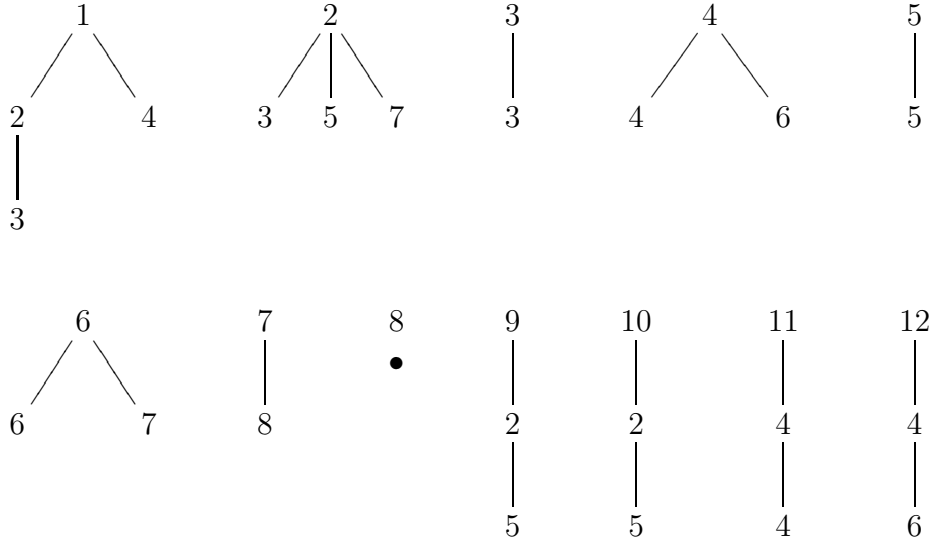
\ignore{
$$\xymatrixcolsep{1pc}
\xy\xymatrix{
 &1  
\ar@{-}[dl] 
\ar@{-}[dr] 
&&&&2\ar@{-}[dl] 
\ar@{-}[d] 
\ar@{-}[dr] 
&&&3\ar@{-}[d] 
&&&4\ar@{-}[dl] 
\ar@{-}[dr] 
&&&5\ar@{-}[d]\\ 
2\ar@{-}[d] 
&&4 &&3&5&7 &&3 &&4&&6 &&5\\
3\\
&6\ar@{-}[dl] \ar@{-}[dr]
&&&7\ar@{-}[d]&& 8\save+<0ex,-3ex>\drop{\bullet}\restore && 9\ar@{-}[d]
&&10\ar@{-}[d]
&&11\ar@{-}[d]
&&12\ar@{-}[d]\\
6&&7 &&8 &&&& 2\ar@{-}[d]&&2\ar@{-}[d]&&4\ar@{-}[d]&&4\ar@{-}[d] \\
&&&& &&&& 5&&5&&4&&6
}\endxy$$
}
 
We will sketch a proof for the fact that $\pinf$ is 
contravariantly finite in $\lamod$ by exhibiting (minimal) 
$\pinf$-approximations of the simple left 
$\Lambda$-modules and by partly verifying that they have the claimed 
properties. Due to \cite{AR}, securing approximations for the simple modules 
guarantees contravariant finiteness. 
The following module $A_1$ is a $\pinf$-approximation of $S_1$. 
Namely, $A_1=B_1\oplus C_1$, where $B_1$ corresponds to the 
representation 

\ignore{
$$\xymatrixcolsep{4pc}
\xy\xymatrix{
K\ar[dr]^(.4){\left({{0} \atop{1}}\right)}\\
K\ar[r]^(.4){\left({{1} \atop {0}}\right)}&K^2\ar@(ur,ul)[]_{0}\ar[r]^\mu&
0\ar@(ur,ul)[]\ar[r]&0\ar[r]&0\\
K\ar[ur]_{\left({{1} \atop {1}}\right)} \ar[dr]^1 \\
K \ar[r]^(0.4)1&K \ar[uurr]^\nu \ar[r] \ar[dr] &0
\ar@(dr,ur)[]\\
0 \ar[ur] &&0 \ar@(dr,ur)[]
}\endxy$$
}

\noindent of $\Gamma$ modulo $I$. The module $C_1$ is defined similarly, with
the  vertex 10 
taking over the role of 9 and the arrow $\gamma$ taking over that of 
$\beta$.

Note first that $\Omega^1(A_1)=(\Lambda e_2/\Lambda\nu)^2 \oplus (\Lambda 
e_4)^2$, hence $\pdim A_1=3$. In particular, $A_1 \in 
\pinf$.
Moreover, denote by $K$ the kernel of the homomorphism $(\psi,\rho) :  
B_1 \oplus C_1 \rightarrow S_1$, where $ \psi :B_1  \rightarrow S_1$ and 
$\rho: C_1 \rightarrow S_1$ are the canonical epimorphisms (unique up to a 
nonzero scalar factor). It can be checked that $\Ext_\Lambda^1(M,K)=0$ for 
each  object $M \in \pinf$, which shows that the map 
 $\Hom_\Lambda (M,A_1) \rightarrow \Hom_\Lambda (M,S_1)$  induced by 
$(\psi,\rho)$ is 
onto, as required. Along the same line, one can verify that the minimal 
$\pinf$-approximations $A_i$ of the simple modules 
$S_i$ with  $i \geq 2$ are as follows: $A_2 =\Lambda e_2/\Lambda \nu$;  
$A_3=\Lambda e_3$; $ A_4$ corresponds to the representation 

\ignore{
$$\xymatrixcolsep{4pc}
\xy\xymatrix{
K\ar[dr]^(.4){\left({{{0} \atop{1}}\atop {0}}\right)}\\
K\ar[r]^(.35){\left({{{0} \atop 
{0}}\atop{1}}\right)}&K^3\ar@(ur,ul)[]_{f_4}\ar[r]^{(1,0,0)}&
K\ar@(ur,ul)[]_0\ar[r]&0\ar[r]&0\\
0\ar[ur] \ar[dr] \\
0 \ar[r]&0 \ar[uurr] \ar[r] \ar[dr] &0
\ar@(dr,ur)[]\\
0 \ar[ur] &&0 \ar@(dr,ur)[]
}\endxy$$
}

\noindent where $f_4$ has matrix 
$\left(\def\objectstyle{\scriptstyle}
\vcenter{\xy\xymatrix@-10mm{0&0&0\\1&0&0\\1&0&0}\endxy}\right)$; moreover, 
$A_5=\Lambda e_5$; $A_6= 
\Lambda e_6/\Lambda \sigma$; $A_7=S_7$; $A_8=S_8$; $A_9$ is $B_1$ as 
defined above; $A_{10}$ is $C_1$ as defined above; $A_{11}=A_{12}=\Lambda 
e_11 \oplus\Lambda e_12/\Lambda(\varepsilon,\chi)$.
Thus Corollary 4.3 is applicable. It yields 
$$\lfindim \Lambda= 
\lFindim \Lambda = \sup \{\pdim A_i \mid 1 \leq i \leq 12\} =\pdim A_1 =3.$$
 
Our final example exhibits a class of finite dimensional monomial 
relation algebras $\Lambda$ such that $\lfindim\Lambda<\lFindim \Lambda$, 
while the subcategory $\pinfrt\subseteq \modla$ of 
the category of all finitely generated right $\Lambda$-modules is 
contravariantly finite.  
\medskip

\noindent {\bf Example 5.2.} Given a positive integer $n$, let $\Lambda=
K\Gamma/I$ be the  monomial relation algebra defined in \cite{HZ:2}, proof of
Theorem E, where 
$\Gamma$ is the quiver

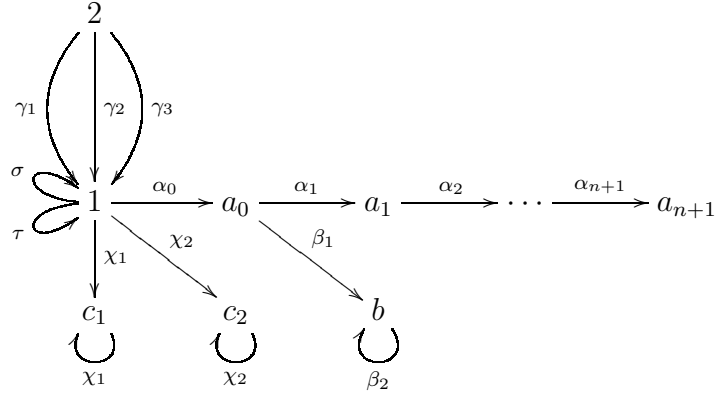
\ignore{
$$\xymatrixcolsep{3pc}
\xy\xymatrix{
2
\ar[dd]^{\gamma_2}
\ar@/^1.5pc/[dd]^{\gamma_3}
\ar@/_1.5pc/[dd]_{\gamma_1}\\
\\
1\ar@(l,ul)[]^{\sigma}
\ar@(l,dl)[]_{\tau}
\ar[d]^{\chi_1}
\ar[dr]^{\chi_2}
\ar[r]^{\alpha_0}&
a_0
\ar[r]^{\alpha_{1}}\ar[dr]^{\beta_{1}}&
a_{1}\ar[r]^{\alpha_{2}}&
\cdots\ar[r]^-{\alpha_{n+1}}&
a_{n+1}\\
c_1\ar@(dr,dl)[]^{\chi_1}&
c_2\ar@(dr,dl)[]^{\chi_2}&
b\ar@(dr,dl)[]^{\beta_2}
}\endxy$$
}

\noindent The relations are such that the indecomposable 
projective right $\Lambda$-modules have the following graphs (in particular,
$e_2\Lambda \cong S_2$).

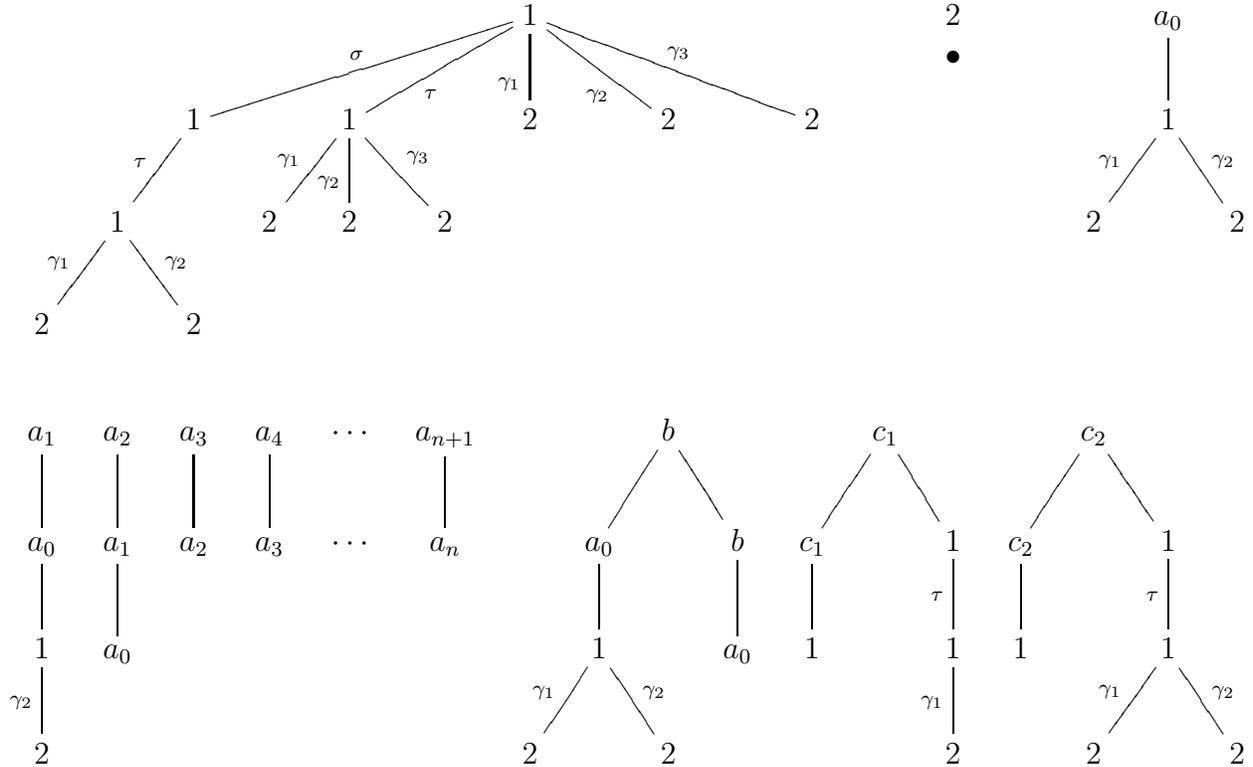
\ignore{ 
$$\xymatrixcolsep{.9pc}
\xy\xymatrix{
&&&&&&1 
\ar@{-}[dllll]_{ \sigma  }
\ar@{-}[dll]^(.6)\tau  
\ar@{-}[d]_(.64){\gamma_1}
\ar@{-}[drr]_(.6){ \gamma_2 }
\ar@{-}[drrrr]^{ \gamma_3}
&&&&&& 2\save+<0ex,-3ex>\drop{\bullet}\restore &&&a_0
\ar@{-}[d]\\
&&1
\ar@{-}[dl]_\tau &&1
\ar@{-}[dl]_{\gamma_1}
\ar@{-}[d]_(.6){\gamma_2}
\ar@{-}[dr]^{\gamma_3} &&2&&2&&2 &&&&&1
\ar@{-}[dl]_{\gamma_1}
\ar@{-}[dr]^{\gamma_2}\\
&1
\ar@{-}[dl]_{\gamma_1}
\ar@{-}[dr]^{\gamma_2} &&2&2&2 &&&&&&&&&2&&2\\
2&&2\\
a_1 \ar@{-}[d] &a_2 \ar@{-}[d] &a_3
\ar@{-}[d]
&a_4
\ar@{-}[d]
&\cdots&a_{n+1}
\ar@{-}[d]
&&&b
\ar@{-}[dl]
\ar@{-}[dr]&&&c_1
\ar@{-}[dl]
\ar@{-}[dr]&&&c_2
\ar@{-}[dl]
\ar@{-}[dr]\\
a_0 \ar@{-}[d] &a_1 \ar@{-}[d] &a_2&a_3&\cdots&a_n&&a_0
\ar@{-}[d]&&b
\ar@{-}[d]&c_1
\ar@{-}[d]&&1
\ar@{-}[d]_\tau&c_2
\ar@{-}[d]&&1
\ar@{-}[d]_\tau\\
1 \ar@{-}[d]_{\gamma_2} &a_0 &&&&&&1
\ar@{-}[dl]_{\gamma_1}
\ar@{-}[dr]^{\gamma_2}&&a_0&1&&1
\ar@{-}[d]_{\gamma_1}&1&&1
\ar@{-}[dl]_{\gamma_1}
\ar@{-}[dr]^{\gamma_2}\\
2&&  &&&&2&&2&&&&2&&2&&2}\endxy$$
}

In view of the fact that there are unique arrows $1 \rightarrow a_0$, $a_i 
\rightarrow a_{i+1}$ for $i \geq 0 $, and $ b \rightarrow b$, $ b 
\rightarrow 
a_0$, $ c_i \rightarrow c_i$, $1 \rightarrow c_i$ for $ i=1,2$, the 
above graphs determine $\Lambda $ up to  isomorphism.

In \cite{HZ:2}, it was shown that $ \lfindim \Lambda = n+1$, 
whereas $\lFindim \Lambda = n+2$. On the other hand, the subcategory 
$\pinfrt$  of finitely generated right $\Lambda$-modules of
finite  projective dimension is contravariantly finite in $\modla$. The 
minimal $\pinfrt$-approximations $A(i)$ of the simple right 
$\Lambda$-modules $e_i\Lambda/e_iJ$ for $i \in \Gamma_0$ are listed 
below. 

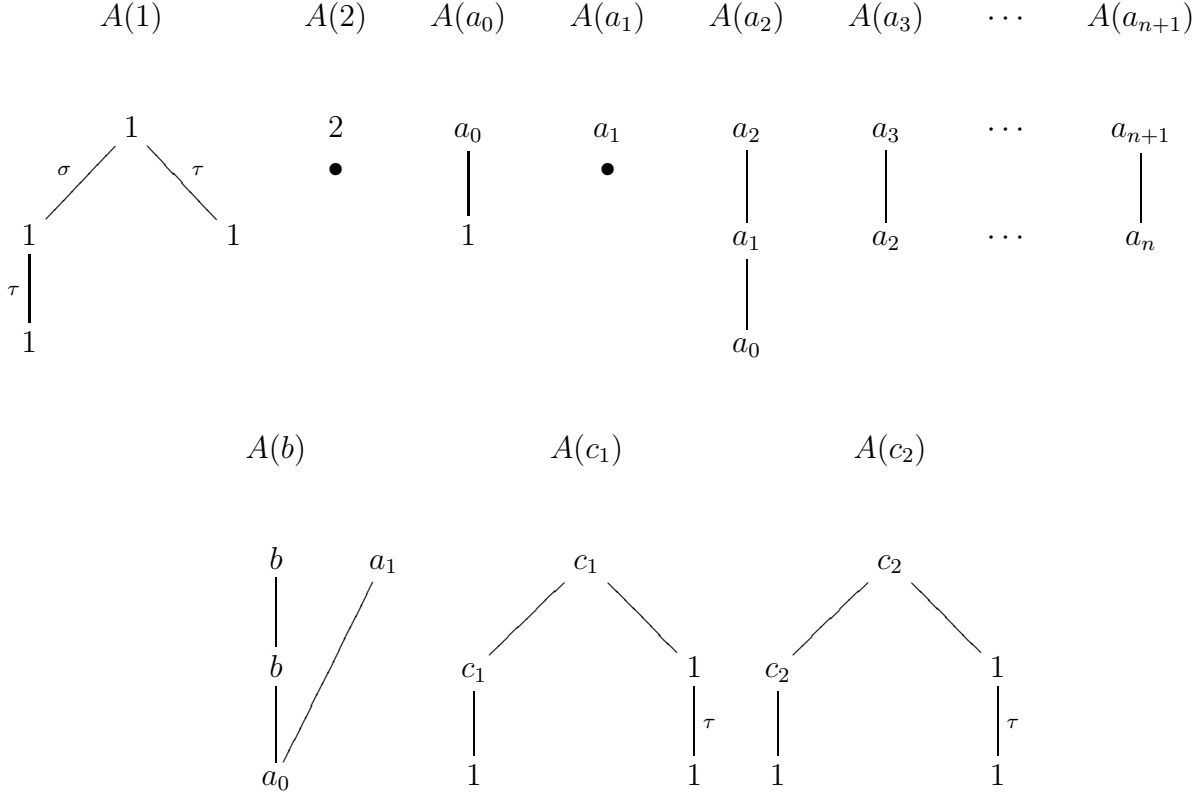
\ignore{
$$\xymatrixcolsep{1.4pc}
\xy\xymatrix{
\ &A(1)&&A(2)&A(a_0)&A(a_1)&A(a_2)&A(a_3)&\cdots&A(a_{n+1})\\
&1\ar@{-}[dl]_\sigma \ar@{-}[dr]^\tau
&&2\save+<0ex,-3ex>\drop{\bullet}\restore  &a_0 \ar@{-}[d] 
&a_1 \save+<0ex,-3ex>\drop{\bullet}\restore &a_2\ar@{-}[d] &a_3
\ar@{-}[d]&\cdots&a_{n+1}
\ar@{-}[d]\\
1
\ar@{-}[d]_\tau&&1&&1&&a_1
\ar@{-}[d]&a_2&\cdots&a_n\\
1&&&&&&a_0&&&\\
}\endxy$$
}

\ignore{
$$\xymatrixcolsep{1.4pc}
\xy\xymatrix{
&A(b)&&&A(c_1)&&&A(c_2)\\
&b
\ar@{-}[d]&a_1
\ar@{-}[ddl]&&c_1
\ar@{-}[dl]
\ar@{-}[dr]&&&c_2
\ar@{-}[dl]
\ar@{-}[dr]\\
&b
\ar@{-}[d]&&c_1
\ar@{-}[d]&&1
\ar@{-}[d]^\tau&c_2
\ar@{-}[d]&&1
\ar@{-}[d]^\tau\\
&a_0&&1&&1&1&&1 
}\endxy$$
}

We give a bit of detail to back our claim for $A(1)$ and $A(b)$. Clearly, 
$\pdim A(1) =1$. Let $M$ be any object in $\pinfrt$ containing a top 
element $x$ of type $e_1$. Since $\pdim e_2\Lambda/e_2J < \infty $ 
and since the vertex $2$ is a source of $\Gamma$, we may assume that 
$Me_2=0$. Note moreover that 
$\Omega^1(M)e_1J^2e_1=0$, because $\Omega^1(M)$ is contained in the 
radical of a projective module. This implies that $\Omega^1(M)$ cannot 
contain a top element of type $e_1$, for otherwise all higher 
syzygies $\Omega^j(M)$ would contain such a top element as well. 
Consequently, the graph of $\Lambda x$ coincides with that of $A(1)$ and 
$\Lambda x $ is a direct summand of $M$. This guarantees that each 
homomorphism $M\rightarrow S_1$ factors through $A(1)$.

We are still briefer in justifying our claim for $A(b)$. We just point to 
the fact that all homomorphisms $M\rightarrow e_b\Lambda/e_bJ$ where 
$M \in \pinfrt$ has a graph belonging to the series

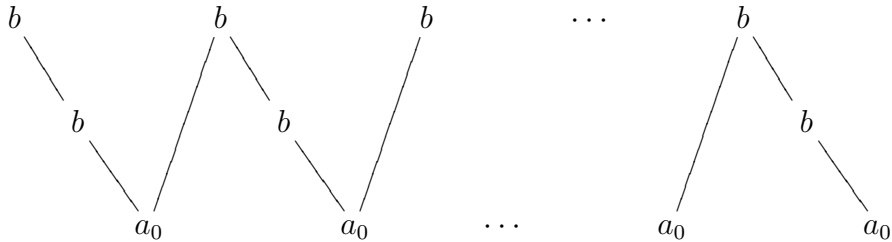
\ignore{
$$\xymatrixcolsep{1pc}
\xy\xymatrix{
b
\ar@{-}[dr]&&&b
\ar@{-}[ddl]
\ar@{-}[dr]
&&&b
\ar@{-}[ddl]
&&\cdots&&b
\ar@{-}[ddl]
\ar@{-}[dr] \\&b
\ar@{-}[dr]&&&b
\ar@{-}[dr]&&&&&&&b
\ar@{-}[dr]\\
&&a_0&&&a_0&&\cdots&&a_0&&&a_0}\endxy$$
}

\noindent factor through $A(b)$. 
Another application of Corollary 4.3 therefore yields that 
$$r \Findim \Lambda= \sup \{\pdim A_i \mid i \in  \Gamma_0 \} =\pdim 
A(a_0)=\pdim A(b)=2.$$

\end{document}